\documentclass[a4paper,12pt]{article}
\usepackage{amsmath}
\usepackage{color}
\usepackage{amssymb}
\usepackage{latexsym}
\usepackage{url}
\numberwithin{equation}{section}

\def\R{{\bf R}}

\def\N{{\bf N}}

\def\d{\displaystyle}
\def\e{{\varepsilon}}
\def\o{\overline}
\def\wt{\widetilde}

\newtheorem{thm}{Theorem}[section]

\newtheorem{lem}{Lemma}[section]

\newtheorem{rem}{Remark}[section]
\newtheorem{Def}{Definition}[section]

\title{Blow-up for semilinear damped wave equations\\
with sub-Strauss exponent\\
in the scattering case
}
\author{
Ning-An Lai
\footnote{Department of Mathematics,
Lishui University,
Lishui City 323000,
China.
e-mail: hyayue@gmail.com.}
\quad
Hiroyuki Takamura
\footnote{Department of Complex and Intelligent Systems,
Faculty of Systems Information Science,
Future University Hakodate,
116-2 Kamedanakano-cho,
Hakodate, Hokkaido 041-8655, Japan.
e-mail: takamura@fun.ac.jp.}
}
\date{
\[
\begin{array}{ll}
\mbox{\footnotesize{\bf Keywords:}}
& \mbox{\footnotesize damped wave equation, semilinear, blow-up, lifespan}\\
\mbox{\footnotesize{\bf MSC2010:}}
& \mbox{\footnotesize primary 35L71, secondary 35B44}\\
\end{array}
\]
}

\begin{document}
\maketitle
\begin{abstract}
It is well-known that
the critical exponent for semilinear damped wave equations
is Fujita exponent when the damping is effective.
Lai, Takamura and Wakasa in 2017 have obtained a blow-up result
not only for super-Fujita exponent  but also for the one closely related to Strauss exponent
when the damping is scaling invariant and its constant is relatively small,
which has been recently extended by Ikeda and Sobajima.
\par
Introducing a multiplier for the time-derivative of the spatial integral of unknown functions,
we succeed in employing the technics on the analysis for  semilinear wave equations
and proving a blow-up result for semilinear damped wave equations
with sub-Strauss exponent when the damping is in the scattering range.
\end{abstract}


\section{Introduction}
\par\quad
In this paper, we consider the following initial value problem
\begin{equation}
\label{1}
\left\{
\begin{array}{l}
\d u_{tt}-\Delta u+\frac{\mu}{(1+t)^\beta}u_t=|u|^p
\quad \mbox{in}\ \R^n\times[0,\infty),\\
u(x,0)=\e f(x),\ u_t(x,0)=\e g(x), \quad  x\in\R^n,
\end{array}
\right.
\end{equation}
where $\mu>0,\ n\in\N,$ $\beta\in\R$ and $\beta>1$.
We assume that $\e>0$ is a \lq\lq small" parameter.
\par
First, we shall outline a background of (\ref{1}) briefly
according to the classifications of Wirth \cite{Wir1, Wir2, Wir3} for the corresponding linear problem.
Let $u^0$ be a solution of the initial value problem for the following linear damped wave equation
\begin{equation}
\label{2}
\left\{
\begin{array}{l}
\d u^0_{tt}-\Delta u^0+\frac{\mu}{(1+t)^\beta}u^0_t=0,
\quad \mbox{in}\ \R^n \times[0,\infty),\\
u^0(x,0)=u_1(x),\ u^0_t(x,0)=u_2(x), \quad x\in\R^n,
\end{array}
\right.
\end{equation}
where $\mu>0$, $\beta\in\R$, $n\in\N$ and $u_1,u_2\in C^\infty_0(\R^n)$.
When $\beta\in(-\infty,-1)$, we say that the damping term is \lq\lq overdamping",
in which case the solution does not decay to zero when $t\rightarrow\infty$.
When $\beta\in[-1,1)$, the solution behaves like that of the heat equation,
which means that the term $u^0_{tt}$ in (\ref{2}) has no influence
on the behavior of the solution.
In fact, $L^p$-$L^q$ decay estimates of the solution
which are almost the same as those of the heat equation are established.
In this case, we say that the damping term is \lq\lq effective."
In contrast, when $\beta\in(1,\infty)$,
it is known that the solution behaves like that of the wave equation,
which means that the damping term in (\ref{2}) has no influence
on the behavior of the solution.
In fact, in this case the solution scatters to that of the free wave equation
when $t\rightarrow\infty$,
and thus we say that we have \lq\lq scattering."
When $\beta=1$, the equation in (\ref{2}) is invariant under the following scaling
\[
\wt{u^0}(x,t):=u^0(\sigma x, \sigma(1+t)-1),\ \sigma>0,
\]
and hence we say that the damping term is \lq\lq scale invariant."
The remarkable fact in this case is that the behavior of the solution of (\ref{2})
is determined by the value of $\mu$.
Actually, for $\mu\in(0,1)$, it is known that the asymptotic behavior of the solution
is closely related to that of the free wave equation.
For this range of $\mu$,
we say that the damping term is \lq\lq non-effective."
However, the threshold of $\mu$ according to the behavior of the solution is still open.
In this way, we may summarize all the classifications of the damping term in (\ref{2})
in the following table.

\begin{center}
\begin{tabular}{|c|c|}
\hline
$\beta\in(-\infty,-1)$ & overdamping\\
\hline
$\beta\in[-1,1)$ & effective\\
\hline
$\beta=1$ &
\begin{tabular}{c}
scaling invariant\\
$\mu\in(0,1)\Rightarrow$ non-effective
\end{tabular}
\\
 \hline
$\beta\in(1,\infty)$ & scattering\\
\hline
\end{tabular}
\end{center}
\par

If $\beta=0$, then we say that the damping term in \eqref{1}
is classical, or of constant coefficient case.
In this case the equation is a good model to describe the wave propagation with the friction,
such as the telegraph equation, the elastic vibration with the damping
proportional to the velocity and the heat conduction with the finite speed of the propagation.
There is extensive literatures on the question of the blow-up in a finite time
or global-in-time existence of the Cauchy problem of semilinear damped wave equation
with constant coefficients.
Based on these works we now know that it admits a critical power,
the so-called Fujita exponent defined by
\begin{equation}
\label{Fujita}
p_F(n):=1+\frac{2}{n},
\end{equation}
which means that the solution will blow up in a finite time for $p\in (1, p_F(n)]$,
and there is a global solution for $p>p_F(n)$ with small data.
We list the related results (but maybe not all of them) in the following table.
\begin{center}
\begin{tabular}{|c|c|c|c|}
\hline
 $1<p<p_F(n)$ & $p=p_F(n)$ & $p>p_F(n)$\\
\hline
\begin{tabular}{l}
Li $\&$ Zhou \cite{LZ}\\($n=1, 2$),\\
Todorova $\&$\\
 Yordanov \cite{TY}
\end{tabular}
 & \begin{tabular}{l}
Li $\&$ Zhou \cite{LZ}\ ($n=1, 2$),\\
Zhang \cite{Z}, or indep.,\\
Kirane $\&$ Qafsaoui \cite{MM}
\end{tabular}  &
\begin{tabular}{l}
Todorova $\&$\\ Yordanov \cite{TY}\\
($p\leq n/(n-2)$\\for $n\ge3$)
\end{tabular}
\\
\hline
\end{tabular}
\end{center}

If the solution blows up in a finite time, people then are interested in the lifespan estimate,
the maximal existence time of the energy solution of (\ref{1}) for arbitrarily fixed $(f,g)$.
We denote it by $T(\e)$.
Now we know the estimate of lifespan will be
\begin{equation}
\label{2a}
T(\e)\sim
\left \{
\begin{array}{cl}
\exp\left(C\varepsilon^{-(p-1)}\right) & \mbox{for}\ p=p_F(n),\\
C\varepsilon^{-2(p-1)/(2-n(p-1))} & \mbox{for}\ 1<p<p_F(n),
\end{array}
\right.
\end{equation}
where $T(\e)\sim A(C,\e)$ stands for the fact that
the following estimate holds with positive constants $C_1,C_2$ independent of $\e$;
\[
A(C_1,\e)\le T(\e)\le A(C_2,\e).
\]
 To our best knowledge, we have the following table,
\begin{center}
{\small
\begin{tabular}{|c|c|c|c|}
\hline
$\beta=0$ & $1<p<p_F(n)$ & $p=p_F(n)$ \\
\hline
$n\le2$ &
\begin{tabular}{l}
U: Li $\&$ Zhou \cite{LZ}, or indep.,\\
\qquad Ikeda $\&$ Wakasugi \cite{IW}\\
L: Fujiwara $\&$ Ikeda $\&$ Wakasugi \cite{FIW}
\end{tabular}
&
\begin{tabular}{l}
U: Li $\&$ Zhou \cite{LZ}\\
L: Ikeda $\&$ Ogawa \cite{IO}
\end{tabular}
\\
\hline
$n=3$ &
\begin{tabular}{l}
U: Ikeda $\&$ Wakasugi \cite{IW}\\
L: Fujiwara $\&$ Ikeda $\&$ Wakasugi \cite{FIW}
\end{tabular}
&
\begin{tabular}{l}
U: Nishihara \cite{N}\\
L: Ikeda $\&$ Ogawa \cite{IO}
\end{tabular}
\\
\hline
$n\ge4$ &
\begin{tabular}{l}
U: Ikeda and Wakasugi \cite{IW}\\
L: Fujiwara $\&$ Ikeda $\&$ Wakasugi \cite{FIW}
\end{tabular}
&
\begin{tabular}{l}
U: Lai $\&$ Zhou \cite{LaZ}\\
L: Ikeda $\&$ Ogawa \cite{IO}
\end{tabular}
\\
\hline
\end{tabular}
}
\end{center}
where \lq\lq U", or \lq\lq L", denotes the upper bound, or the lower bound, of the lifespan estimate respectively.

Recently, the Cauchy problem for the damped wave equation with variable coefficients
and the power type source term, \eqref{1},
attracts more and more attention. We want to see whether this equation still has a diffusive structure. Generally speaking, based on the classification for the linear problem, \eqref{2}, mentioned above, one may expect that if the coefficient $\mu/(1+t)^\beta$ decays not so fast, then the damping is effective, i.e. the solution behaves like that of the corresponding nonlinear heat equation. And if $\mu/(1+t)^\beta$ decays sufficiently fast, then the damping becomes non-effective, i.e. the solution behaves like that of the nonlinear wave equation. Essentially, it is totally determined by the constants $\mu$ and $\beta$. Till now, most of the known results focus on the range $\beta\in (-\infty, 1]$.
\par
 When $\beta\in (-1, 1)$, the global existence, i.e. $T(\e)=\infty$, has been obtained by
 D'Abbicco, Lucente and Reissig \cite{DLR13} for  $p_F(n)<p<n/[n-2]_+$, where
 \[
 \frac{m}{[n-2]_+}:=
\left\{ \begin{array}{cl}
 \infty & \mbox{for}\ n=1,2,\\
 m/(n-2) & \mbox{for}\ n\ge3.
 \end{array}
 \right.
 \]
Nishihara \cite{N1} has extended this upper bound to $(n+2)/[n-2]_+$ for $\beta\in[0,1)$.
Lin, Nishihara and Zhai \cite{LNZ12} have removed this restriction on $\beta$.
When $\beta=-1$, Wakasugi \cite{WY17} has obtained the global existence
for $p_F(n)<p<n/[n-2]_+$.
In the counter parts on $\beta\in[-1,1)$,
we have the following lifespan estimates by
Ikeda and Ogawa \cite{IO}, Fujiwara, Ikeda and Wakasugi \cite{FIW},
Ikeda and Inui \cite{II}.
\begin{center}
\begin{tabular}{|c|c|c|}
\hline
$T(\e)$ & $1<p<p_F(n)$ & $p=p_F(n)$\\
\hline
$\beta=-1$ &
 $\sim\exp\left(C\e^{-2(p-1)/(2-n(p-1))}\right)$ \cite{FIW}
&
\begin{tabular}{c}
$\sim\exp\left(\exp(C\e^{-(p-1)})\right)$\\
L: \cite{FIW}, U: \cite{II}
\end{tabular}
\\
\hline
$
\begin{array}{c}
\beta\in(-1, 1)\\
(\beta\neq0)
\end{array}
$
&
$\sim C\e^{-2(p-1)/\{(2-n(p-1))(1+\beta)\}}$ \cite{FIW}
&
\begin{tabular}{c}
 $\sim\exp\left(C\varepsilon^{-(p-1)}\right)$\\
L: \cite{FIW}, U: \cite{II}
\end{tabular}
\\
\hline
\end{tabular}
\end{center}
So we can conclude that the critical exponent is still $p_F(n)$ in (\ref{Fujita})
which is the one for semilinear heat equation, $u_t-\Delta u=|u|^p$.
When $\beta<-1$, Ikeda and Wakasugi \cite{IWnew} have recently proved that
the global existence holds for any $p>1$.
\par
 When $\beta=1$, the situation is a bit complicated.
 It seems to be the threshold between effective and non-effective damping.
 Actually, the other constant $\mu$ also plays a crucial role in this case.
 D'Abbicco and Lucente \cite{DL1}, and D'Abbicco \cite{DABI} have showed that
 the critical power is $p_F(n)$ when
 \[
 \mu\geq
 \left\{
 \begin{array}{cl}
5/3 & \mbox{for}\ n=1,\\
3 & \mbox{for}\ n=2,\\
n+2 & \mbox{for}\ n\geq 3.
\end{array}
\right.
\]
Noting that $\mu=2$ is an exceptional case, since after making
the so-called Liouville transform
\begin{equation}\nonumber
\begin{aligned}
w(x,t):=(1+t)^{\mu/2}u(x,t),
\end{aligned}
\end{equation}
then problem (\ref{1}) can be rewritten as
\begin{equation}
\label{3}
\left\{
\begin{array}{l}
\d w_{tt}-\Delta w+\frac{\mu(2-\mu)}{4(1+t)^2}w=\frac{|w|^p}{(1+t)^{\mu(p-1)/2}}
\quad \mbox{in}\ \R^n \times[0,\infty),\\
w(x,0)=\e f(x),\ w_t(x,0)=\e \{(\mu/2)f(x)+g(x)\}, \quad x\in\R^n.
\end{array}
\right.
\end{equation}
When $\beta=1$ and $\mu=2$,
it is natural to think that the critical power is related to the so-called Strauss exponent $p_S(n)$ which is defined for $n>1$ as a positive root of the quadratic equation,
\begin{equation}
\label{4}
\gamma(p,n):=2+(n+1)p-(n-1)p^2=0.
\end{equation}
We note that
\[
p_F(n)<p_S(n)=\frac{n+1+\sqrt{n^2+10n-7}}{2(n-1)}
\quad\mbox{for}\ n\ge2
\]
and that $p_S(n)$ is the critical power for semilinear wave equations,
$u_{tt}-\Delta u=|u|^p$.
D'Abbicco, Lucente and Reissig \cite{DLR15} have determined
a critical power
\[
p_c(n):=\max\{p_F(n), p_S(n+2)\}\ \mbox{for}\ n\leq 3.
\]
D'Abbicco and Lucente \cite{D-L} have also showed the global existence for $p_S(n+2)<p<1+2/\max\{2, (n-3)/2\}$ to odd and higher dimensions($n\geq 5$)
under the spherically symmetric assumption.
In the blow-up case, we have results about the lifespan estimate
as in the following table.
\begin{center}
\begin{tabular}{|c||c|c|}
\hline
$T(\e)$ & $1<p<p_c(n)$ & $p=p_c(n)$\\
\hline
\hline
$n=1$
&
$
\sim
\left\{
\begin{array}{l}
C\e^{-(p-1)/(3-p)}\\
\qquad\mbox{for}\ I_{f,g}\neq0,\\
C\e^{-p(p-1)/(3-p)}\\
\qquad\mbox{for $\ I_{f,g}=0$ and $2<p$},\\
Cb(\e) \\
\qquad\mbox{for $I_{f,g}=0$ and $p=2$},\\
C\e^{-2p(p-1)/\gamma(p,3)}\\
\qquad\mbox{for $I_{f,g}=0$ and $p<2$}
\end{array}
\right.$
&
$\sim
\left\{
\begin{array}{l}
\exp\left(C\e^{-(p-1)}\right)\\
\qquad\mbox{for}\  I_{f,g}\neq0,\\
\exp\left(C\e^{-p(p-1)}\right)\\
\qquad\mbox{for}\  I_{f,g}=0
\end{array}
\right.
$
\\
\hline
$n=2$
&
$\sim
\left\{
\begin{array}{l}
C\e^{-(p-1)/(4-2p)}\\
\qquad\mbox{for}\  I_{f,g}\neq0,\\
C\e^{-2p(p-1)/\gamma(p,4)}\\
\qquad\mbox{for}\  I_{f,g}=0
\end{array}
\right.
$
&
$\sim
\left\{
\begin{array}{l}
\exp\left(C\e^{-(p-1)}\right)\\
\qquad\mbox{for}\  I_{f,g}\neq0,\\
\exp\left(C\e^{-p(p-1)}\right)\\
\qquad\mbox{for}\ I_{f,g}=0
\end{array}
\right.
$
\\
\hline
$n=3$ & $\sim C\e^{-2p(p-1)/\gamma(p,n+2)}$
& $\sim  \exp\left(C\e^{-p(p-1)}\right)$
\\
\hline
\end{tabular}
\end{center}
where $b=b(\e)$ is a positive number satisfying $b\e^2\log^2(b+1)=1$ and
\[
I_{f,g}:=\int_{\R^n}\{f(x)+g(x)\}dx.
\]
See Wakasa \cite{wak16}, and Kato, Takamura, Wakasa \cite{KTW} for $n=1$,
Imai, Kato, Takamura and Wakasa \cite{IKTW} for $n=2$,
Kato and Sakuraba \cite{KS} for $n=3$.
\par
When $\beta=1$ and $\mu\neq2$,
due to the work of Wakasugi \cite{WY14_scale}, we may believe that the solution is \lq\lq heat-like"(the critical power is Fujita exponent) for $\mu>1$.
Actually he established the following lifespan estimates,
\begin{equation}
\label{lifespan:wakasugi}
\left\{
\begin{array}{l}
T(\e)\le C\e^{-(p-1)/\{2-n(p-1)\}}\\
\qquad\mbox{for $1<p<p_F(n)$ and $\mu\ge1$},\\
T(\e)\le C\e^{-(p-1)/\{2-(n+\mu-1)(p-1)\}}\\
\qquad\mbox{for $1<p<p_F(n+\mu-1)$ and $0<\mu<1$}
\end{array}
\right.
\end{equation}
as well as $T(\e)<\infty$ for $p=p_F(n)$ and $\mu\ge1$,
and for $p=p_F(n+\mu-1)$ and $0<\mu<1$.
Surprisingly, the authors and Wakasa \cite{LTW} found that the solution is \lq\lq wave-like"(the critical exponent is bigger than Fujita exponent
and is related to the Strauss exponent) in some case even for $\mu>1$, by using the improved Kato's lemma introduced by the second author \cite{Takamura15}.
Actually we have the estimate,
\begin{equation}
\label{lifespan:LTW}
\begin{array}{l}
T(\e)\le C\e^{-2p(p-1)/\gamma(p,n+2\mu)}\\
\qquad\mbox{for $p<p_S(n+2\mu)$ and
$\d0<\mu<\frac{n^2+n+2}{2(n+2)}$}.
\end{array}
\end{equation}
We note that (\ref{lifespan:LTW}) is stronger than (\ref{lifespan:wakasugi})
for $n\ge2$.
Very recently, we have been informed that Ikeda and Sobajima \cite{IS} extend this result to
\[
T(\e)\le
\left\{
\begin{array}{ll}
\exp\left(C\e^{-p(p-1)}\right) & \mbox{for}\ p=p_S(n+\mu),\\
C_\delta\e^{-2p(p-1)/\gamma(p,n+\mu)-\delta} &
\mbox{for}\ p_S(n+2+\mu)\le p<p_S(n+\mu),\\
C_\delta\e^{-1-\delta} & \mbox{for}\ p_F(n)<p<p_S(n+2+\mu)
\end{array}
\right.
\]
when
\[
n\ge2\quad\mbox{and}\quad0\le\mu<\frac{n^2+n+2}{n+2},
\]
and
\[
T(\e)\le
\left\{
\begin{array}{ll}
\exp\left(C\e^{-p(p-1)}\right) & \mbox{for}\ p=p_S(1+\mu),\\
C_\delta\e^{-2p(p-1)/\gamma(p,1+\mu)-\delta} &
\mbox{for}\ \max\{3,2/\mu\}\le p<p_S(1+\mu),\\
C_\delta\e^{-2(p-1)/\mu-\delta} & \mbox{for}\ 0<\mu<2/3,\ 3\le p<2/\mu
\end{array}
\right.
\]
when $n=1$ and $0<\mu<4/3$,
with arbitrary small $\delta>0$.

\begin{rem}
For the scale invariant case, $\beta=1$, it is still open to determine the critical exponent. we also remark that for $\beta=1$, the case of $\mu\approx 1$
is related to the semilinear generalized Tricomi equation, which comes from the gas dynamics, see \cite{HWY1, HWY2}.
\end{rem}

In this paper, we are devoted to studying the Cauchy problem \eqref{1} with $\beta>1$. Due to the authors' best knowledge, this problem is completely open. As mentioned above, the corresponding linear problem belongs to the scattering case. We then expect that the solution behaves like that of the semilinear wave equation without the damping term. We are mainly concerned about the blow-up result and upper bound of the lifespan estimate. The novelty is that we introduce a multiplier of exponential type, which is bounded from above and below. Then we get the lower bound of the time-derivative of the spatial integral of the unknown function by the space-time integral of the nonlinear term. Finally, the desired blow-up result and lifespan estimate for sub-Strauss exponent are established by an iteration argument.

\begin{rem}
Compared to the scale invariant case, $\beta=1$, the main difficulty is that we cannot use Liouville transform to rewrite the equation in a form of nonlinear wave or Klein-Gordon equation. See (\ref{3}). We overcome this obstacle by introducing a \lq\lq good" multiplier.
\end{rem}

We organize this paper in five sections.
In Section 2, we give main results.
In Section 3 the key multiplier is introduced,
and the lower bound of the nonlinear term is obtained.
In Section 4, we obtain the blow-up result  
and the upper bound of the lifespan for sub-Strauss exponent
by an iteration argument.
We also give improvements of estimates of the lifespan
for one dimensional case and two dimensional case with low powers in Section 5
under an additional assumption on the initial speed.

\section{Main Result}
\par\quad
As in the work \cite{LTW}, we first define the energy and weak  solution of the Cauchy problem \eqref{1}.
\begin{Def}\label{def1}
We say that $u$ is an energy solution of \eqref{1} over $[0,T)$
if
\begin{equation}
\label{5}
u\in C([0,T),H^1(\R^n))\cap C^1([0,T),L^2(\R^n))\cap L_{\rm loc}
^p(\R^n\times(0,T))
\end{equation}
satisfies
\begin{equation}
\label{6}
\begin{array}{l}
\d\int_{\R^n}u_t(x,t)\phi(x,t)dx-\int_{\R^n}u_t(x,0)\phi(x,0)dx\\
\d+\int_0^tds\int_{\R^n}\left\{-u_t(x,s)\phi_t(x,s)+\nabla u(x,s)\cdot\nabla\phi(x,s)\right\}dx\\
\d+\int_0^tds\int_{\R^n}\frac{\mu u_t(x,s)}{(1+s)^{\beta}}\phi(x,s)dx
=\int_0^tds\int_{\R^n}|u(x,s)|^p\phi(x,s)dx
\end{array}
\end{equation}
with any $\phi\in C_0^{\infty}(\R^n\times[0,T))$ and any $t\in[0,T)$.
\end{Def}

Employing the integration by parts in \eqref{6}
and letting $t\rightarrow T$, we have that
\[
\begin{array}{l}
\d\int_{\R^n\times[0,T)}
u(x,s)\left\{\phi_{tt}(x,s)-\Delta \phi(x,s)
-\left(\frac{\mu\phi(x,s)}{(1+s)^{\beta}}\right)_s \right\}dxds\\
\d=\int_{\R^n}\mu u(x,0)\phi(x,0)dx-\int_{\R^n}u(x,0)\phi_t(x,0)dx\\
\d\quad+\int_{\R^n}u_t(x,0)\phi(x,0)dx+\int_{\R^n\times[0,T)}|u(x,s)|^p\phi(x,s)dxds.
\end{array}
\]
This is exactly the definition of the weak solution of \eqref{1}.

Our main results are stated in the following three theorems.

\begin{thm}
\label{thm1}
Let $\beta>1$ and
\[
 1<p<
 \left\{
 \begin{array}{ll}
 p_S(n) & \mbox{for}\ n\ge2,\\
 \infty & \mbox{for}\ n=1.
 \end{array}
 \right.
 \]
Assume that both $f\in H^1(\R^n)$ and $g\in L^2(\R^n)$ are non-negative,
and $f$ does not vanish identically.
Suppose that an energy solution $u$ of (\ref{1}) satisfies
\begin{equation}
\label{support}
\mbox{\rm supp}\ u\ \subset\{(x,t)\in\R^n\times[0,T)\ :\ |x|\le t+R\}
\end{equation}
with some $R\ge1$.
Then, there exists a constant $\e_0=\e_0(f,g,n,p,\mu, \beta, R)>0$
such that $T$ has to satisfy
\begin{equation}
\label{7}
\begin{aligned}
T\leq C\varepsilon^{-2p(p-1)/\gamma(p, n)}
\end{aligned}
\end{equation}
for $0<\e\le\e_0$, where $C$ is a positive constant independent of $\e$.
\end{thm}

\begin{rem}
In (\ref{7}) for $n=1$, we note that $\gamma(p,1)=2+2p$ by its definition (\ref{4}).
\end{rem}

In low dimensions, we have improvements on the lifespan estimates
as follows.

\begin{thm}
\label{thm2}
Let $n=2$ and $1<p<2$.
Assume that the initial data satisfy the same condition as that in Theorem \ref{thm1}. If
\begin{equation}
\label{7a}
\begin{aligned}
\int_{\R^2}g(x)dx\neq0,
\end{aligned}
\end{equation}
then (\ref{7}) is replaced by
\begin{equation}
\label{7b}
\begin{aligned}
T\leq C\varepsilon^{-(p-1)/(3-p)}.
\end{aligned}
\end{equation}
\end{thm}

\begin{rem}
(\ref{7b}) is stronger than (\ref{7}) by the fact that
$1<p<2$ is equivalent to
\[
\frac{p-1}{3-p}<\frac{2p(p-1)}{\gamma(p,2)}.
\]
\end{rem}

\begin{thm}\label{thm3}
Let $n=1$ and $p>1$. 
Assume that the initial data satisfy the same condition as that in Theorem \ref{thm1}. If
\begin{equation}
\label{7c}
\begin{aligned}
\int_{\R}g(x)dx\neq0,
\end{aligned}
\end{equation}
then  (\ref{7}) is replaced by
\begin{equation}
\label{7d}
\begin{aligned}
T\leq C\varepsilon^{-(p-1)/2}.
\end{aligned}
\end{equation}
\end{thm}

\begin{rem}
(\ref{7d}) is stronger than (\ref{7}) by the fact that
$p>1$ is equivalent to
\[
\frac{p-1}{2}<\frac{2p(p-1)}{\gamma(p,1)}.
\]
\end{rem}

\begin{rem}
These results in the theorems above are the same
as those of Cauchy problem of semilinear wave equations,
$u_{tt}-\Delta u=|u|^p$, except for the case of
\[
n=p=2\quad\mbox{and}\quad\int_{\R^2}g(x)dx\neq0.
\]
\end{rem}
For example, see the introductions in
Takamura \cite{Takamura15} and
Imai, Kato, Takamura and Wakasa \cite{IKTW17}.

\section{Lower Bound of the Functional}
\par\quad
In this section we introduce the multiplier for our problem. Let
\begin{equation}
\label{8}
\begin{aligned}
m(t):=\exp\left(\mu\frac{(1+t)^{1-\beta}}{1-\beta}\right).
\end{aligned}
\end{equation}
Then it is easy to see that for $\beta>1$ we have
\begin{equation}\label{9}
\begin{aligned}
1\ge m(t) \ge m(0) \quad\mbox{for}\ t\ge0,
\end{aligned}
\end{equation}
which means that $m(t)$ is bounded from both above and below.

\begin{rem}
If one puts
\[
m(t)=\exp\left(\mu\log(1+t)\right)=(1+t)^\mu
\]
instead of (\ref{8}) in the proof below,
it will give us a simple proof of the result of
Lai, Takamura and Wakasa \cite{LTW}
which is cited in (\ref{lifespan:LTW}).
\end{rem}

\par
Set
\[
F_0(t):=\int_{\R^n}u(x,t)dx.
\]
Choosing the test function $\phi=\phi(x,s)$ in (\ref{6}) to satisfy
$\phi\equiv 1$ in $\{(x,s)\in \R^n\times[0,t]:|x|\le s+R\}$, we get
\[
\begin{array}{l}
\d\int_{\R^n}u_t(x,t)dx-\int_{\R^n}u_t(x,0)dx\\
\d+\int_0^tds\int_{\R^n}\frac{\mu u_t(x,s)}{(1+s)^\beta}dx
=\int_0^tds\int_{\R^n}|u(x,s)|^pdx,
\end{array}
\]
which means that
\begin{equation}
\label{10}
F_0''+\frac{\mu}{(1+t)^\beta}F_0'
=\int_{\R^n}|u(x,t)|^pdx.
\end{equation}
Multiplying the both sides of \eqref{10} with $m(t)$ yields
\begin{equation}
\label{11}
\left\{m(t)F'_0\right\}'
=m(t)\int_{\R^n}|u(x,t)|^pdx.
\end{equation}
We then get the lower bound of $F'_0(t)$ by integrating \eqref{11} over $[0, t]$
\begin{equation}
\label{12}
F'_0(t)\ge m(0)\int_0^tds\int_{\R^n}|u(x,s)|^pdx
\quad\mbox{for}\ t\ge0,
\end{equation}
where we used the fact that $F_0(0)>0, F'_0(0)\geq0$ and \eqref{9}.

\par\quad
Let
\begin{equation}\nonumber
F_1(t):=\int_{\R^n}u(x,t)\psi_1(x,t)dx,
\end{equation}
where
\begin{equation}\nonumber
\psi_1(x,t):=e^{-t}\phi_1(x),
\quad
\phi_1(x):=
\left\{
\begin{array}{ll}
\d\int_{S^{n-1}}e^{x\cdot\omega}dS_\omega & \mbox{for}\ n\ge2,\\
e^x+e^{-x} & \mbox{for}\ n=1.
\end{array}
\right.
\end{equation}
This is a test function introduced by Yordanov and Zhang \cite{YZ06}.

\begin{lem}[Inequality (2.5) of Yordanov and Zhang \cite{YZ06}]\label{lem1}
\begin{equation}
\label{14}
\int_{|x|\leq t+R}\left[\psi_1(x,t)\right]^{p/(p-1)}dx
\leq C(1+t)^{(n-1)\{1-p/(2(p-1))\}},
\end{equation}
where $C=C(n,p,R)>0$.
\end{lem}

By H\"{o}lder's inequality and \eqref{14}, one has
\begin{equation}\label{15}
\begin{aligned}
\int_{\R^n}|u(x,t)|^pdx\ge C_1(1+t)^{(n-1)(1-p/2)}|F_1(t)|^p
\quad\mbox{for}\ t\ge0,
\end{aligned}
\end{equation}
where $C_1=C_1(n,p,R)>0$. Now we are left with the lower bound of $F_1(t)$. We start with the definition of the energy solution \eqref{6}, which yields that
\[
\begin{array}{l}
\d\frac{d}{dt}\int_{\R^n}u_t(x,t)\phi(x,t)dx\\
\d+\int_{\R^n}\left\{-u_t(x,t)\phi_t(x,t)-u(x,t)\Delta\phi(x,t)\right\}dx\\
\d+\int_{\R^n}\frac{\mu u_t(x,t)}{(1+t)^\beta}\phi(x,t)dx
=\int_{\R^n}|u(x,t)|^p\phi(x,t)dx.
\end{array}
\]
Multiplying the both sides of the above equality with $m(t)$, we have that
\[
\begin{array}{l}
\d\frac{d}{dt}\left\{m(t)
\int_{\R^n}u_t(x,t)\phi(x,t)dx\right\}\\
\d+m(t)\int_{\R^n}\left\{-u_t(x,t)\phi_t(x,t)-u(x,t)\Delta\phi(x,t)\right\}dx\\
\d=m(t)\int_{\R^n}|u(x,t)|^p\phi(x,t)dx.
\end{array}
\]
Integrating this equality over $[0,t]$, we get
\[
\begin{aligned}
& m(t)
\int_{\R^n}u_t(x,t)\phi(x,t)dx
-m(0)\e\int_{\R^n}g(x)\phi(x,0)dx\\
&-\int_0^tds\int_{\R^n}m(s)u_t(x,s)\phi_t(x,s)dx\\
&=\int_0^tds\int_{\R^n}\left\{m(s)
u(x,s)\Delta\phi(x,s)+m(s)|u(x,s)|^p\phi(x,s)\right\}dx.
\end{aligned}
\]

If we put
\[
\phi(x,t)=\psi_1(x,t)=e^{-t}\phi_1(x)
\quad\mbox{on}\quad \mbox{supp}\ u,
\]
we have
\[
\phi_t=-\phi,\ \phi_{tt}=\Delta\phi \quad\mbox{on}\quad \mbox{supp}\ u.
\]
Hence we obtain that
\[
\begin{array}{l}
\d m(t)\{F_1'(t)+2F_1(t)\}-m(0)\e\int_{\R^n}\left\{f(x)+g(x)\right\}\phi(x)dx\\
\d=\int_0^tm(s)
\frac{\mu}{(1+s)^\beta}F_1(s)ds+\int_0^tds\int_{\R^n}m(s)|u(x,s)|^pdx,
\end{array}
\]
which yields
\[
\begin{aligned}
F'_1(t)+2F_1(t)
&\d \ge\frac{m(0)}{m(t)}
C_{f,g}\e+\frac{1}{m(t)}\int_0^tm(s)
\frac{\mu}{(1+s)^\beta}F_1(s)ds\\
&\d \ge m(0)C_{f,g}\e+\frac{1}{m(t)}\int_0^tm(s)
\frac{\mu}{(1+s)^\beta}F_1(s)ds,
\end{aligned}
\]
where
\[
C_{f,g}:=\int_{\R^n}\left\{f(x)+g(x)\right\}\phi_1(x)dx>0.
\]
Integrating the above inequality over $[0,t]$ with a multiplication by $e^{2t}$, we get
\begin{equation}
\label{17}
\begin{array}{ll}
e^{2t}F_1(t)
&\d\ge F_1(0)+m(0)C_{f,g}\e\int_0^te^{2s}ds\\
&\d+\int_0^t\frac{e^{2s}}{m(s)}ds
\int_0^sm(r)\frac{\mu}{(1+r)^\beta}F_1(r)dr.
\end{array}
\end{equation}
Due to a comparison argument, we have that $F_1(t)>0$ for $t\ge0$.
Actually, $F_1(0)>0$ and the continuity of $F_1(t)$ in $t$ yield that
$F_1(t)>0$ for small $t>0$.
If there is the nearest zero point $t_0$ to $t=0$ of $F_1$,
then (\ref{17}) gives a contradiction at $t_0$.

Therefore we obtain that
\[
\begin{array}{ll}
e^{2t}F_1(t)
&\d \ge F_1(0)+m(0)C_{f,g}\e\int_0^te^{2s}ds\\
&\d \ge m(0)F_1(0)
+\frac{1}{2}m(0)C_{f,0}\e(e^{2t}-1)\\
&\d > \frac{1}{2}m(0)C_{f,0}\e e^{2t}
\end{array}
\]
because of
\[
F_1(0)=C_{f,0}\e,\quad C_{f,g}\geq C_{f,0},
\]
which in turn gives us the lower bound of $F_1(t)$,
\begin{equation}
\label{18}
F_1(t)>\frac{1}{2}m(0)C_{f,0}\e
\quad\mbox{for}\ t\ge0.
\end{equation}
\section{Proof for Theorem \ref{thm1}}
\par\quad
By H\"{o}lder inequality with (\ref{support}), it is easy to get
\begin{equation}
\label{19a}
\int_{\R^n}|u(x,t)|^pdx\ge C_2(1+t)^{-n(p-1)}|F_0(t)|^p
\quad\mbox{for}\ t\ge0,
\end{equation}
where $C_2=C_2(n,p,R)>0$.
Then it follows from \eqref{12} and \eqref{19a} that
\begin{equation}
\label{19b}
F_0(t)>C_3\int_0^tds
\int_0^s(1+r)^{-n(p-1)}F_0(r)^pdr
\quad\mbox{for}\ t\ge0,
\end{equation}
where
\[
C_3:=C_2m(0)>0.
\]
Plugging \eqref{18} into \eqref{12} with \eqref{15}, we have
\begin{equation}\label{19}
F'_0(t)\ge C_4\e^p\int_0^t(1+s)^{(n-1)(1-p/2)}ds
\quad\mbox{for}\ t\ge0,
\end{equation}
where
\[
C_4:=C_1m(0)\left(\frac{1}{2}m(0)C_{f,0}\right)^p.
\]
Integrating \eqref{19} over $[0, t]$, we have
\begin{equation}
\label{20}
\begin{array}{ll}
F_0(t)
&\d >C_4\e^p\int_0^tdr
\int_0^r(1+s)^{(n-1)(1-p/2)}ds\\
&\d \ge C_4\e^p(1+t)^{-(n-1)p/2}\int_0^tdr
\int_0^rs^{n-1}ds\\
&\d=\frac{C_4}{n(n+1)}\e^p(1+t)^{-(n-1)p/2}t^{n+1}
\quad\mbox{for}\ t\ge0.
\end{array}
\end{equation}

Next we will begin our iteration argument. First we may assume that $F_0(t)$ satisfies
\begin{equation}
\label{21}
F_0(t)>D_j(1+t)^{-a_j}t^{b_j}\quad\mbox{for}\ t\ge0
\quad(j=1,2,3\cdots)
\end{equation}
with positive constants, $D_j,a_j,b_j$,
which will be determined later.
Due to \eqref{20}, noting that \eqref{21} is true with
\begin{equation}
\label{22}
D_1=\frac{C_4}{n(n+1)}\e^p,
\quad a_1=(n-1)\frac{p}{2},
\quad b_1=n+1.
\end{equation}
Plugging \eqref{21} into \eqref{19b}, we have
\[
\begin{array}{ll}
F_0(t)
&\d>C_3D_j^p\int_0^tds
\int_0^s(1+r)^{-n(p-1)-pa_j}r^{pb_j}dr\\
&\d>C_3D_j^p(1+t)^{-n(p-1)-pa_j}\int_0^tds
\int_0^sr^{pb_j}dr\\
&\d>\frac{C_3D_j^p}{(pb_j+2)^2}(1+t)^{-n(p-1)-pa_j}t^{pb_j+2}
\quad\mbox{for}\ t\ge0.
\end{array}
\]
So we can define the sequences $\{D_j\}$, $\{a_j\}$, $\{b_j\}$ by
\begin{equation}
\label{23}
D_{j+1}\ge\frac{C_3D_j^p}{(pb_j+2)^2},
\quad
a_{j+1}=pa_j+n(p-1),
\quad
b_{j+1}=pb_j+2
\end{equation}
to establish
\[
F_0(t)>D_{j+1}(1+t)^{-a_{j+1}}t^{b_{j+1}}\quad\mbox{for}\ t\ge0.
\]

It follows from \eqref{22} and \eqref{23} that
\[
a_j=p^{j-1}\left((n-1)\frac{p}{2}+n\right)-n
\qquad\mbox{for}\ j=1,2,3,\cdots
\]
and
\[
b_j=p^{j-1}\left(n+1+\frac{2}{p-1}\right)-\frac{2}{p-1}
\qquad\mbox{for}\ j=1,2,3,\cdots.
\]
If we employ the inequality
\[
b_{j+1}=pb_j+2\le p^j\left(n+1+\frac{2}{p-1}\right)
\]
in \eqref{23}, we have
\begin{equation}
\label{24a}
D_{j+1}\ge C_5\frac{D_j^p}{p^{2j}},
\end{equation}
where
\[
C_5:=\frac{C_3}{\left(n+1+\d\frac{2}{p-1}\right)^2}.
\]
From \eqref{24a} it holds that
\[
\begin{aligned}
\log D_j&\geq p\log D_{j-1}-2(j-1)\log p+\log C_5\\
&\geq p^2\log D_{j-2}-2\big(p(j-2)+(j-1)\big)\log p+(p+1)\log C_5.
\end{aligned}
\]
Repeating this procedure, we have
\[
\log D_j
\geq p^{j-1}\log D_1-\sum_{k=1}^{j-1}\frac{2k\log p-\log C_5}{p^k},
\]
which yields that
\[
D_j\ge\exp\left\{p^{j-1}\left(\log D_1-S_p(j)\right)\right\},
\]
where
\[
S_p(j):=\sum_{k=1}^{j-1}\frac{2k\log p-\log C_5}{p^k}.
\]
By d'Alembert's criterion we know that $S_p(j)$ converges for $p>1$ as $j\rightarrow\infty$.
Hence we obtain that
\[
D_j\ge\exp\left\{p^{j-1}\left(\log D_1-S_p(\infty)\right)\right\}.
\]
Turning back to \eqref{21}, we have
\begin{equation}\label{27}
\begin{aligned}
F_0(t)\ge(1+t)^nt^{-2/(p-1)}\exp\left(p^{j-1}J(t)\right)
\quad\mbox{for}\ t>0,
\end{aligned}
\end{equation}
where
\[
\begin{aligned}
J(t)&=-\left((n-1)\frac{p}{2}+n\right)\log(1+t)
+\left(n+1+\frac{2}{p-1}\right)\log t\\
&+\log D_1-S_p(\infty).
\end{aligned}
\]
For $t\geq 1$, by the definition of $J(t)$ we have
\[
\begin{aligned}
J(t)&\geq -\Big((n-1)\frac{p}{2}+n\Big)\log(2t)+\big(n+1+\frac{2}{p-1}\big)\log t\\
&+\log D_1-S_p(\infty)\\
&=\frac{\gamma(p, n)}{2(p-1)}\log t+\log D_1-\Big((n-1)\frac{p}{2}+n\Big)\log 2-S_p(\infty)\\
&=\log \big(t^{\gamma(p, n)/\{2(p-1)\}}D_1\big)-C_6,\\
\end{aligned}
\]
where
\[
 C_6:=\Big((n-1)\frac{p}{2}+n\Big)\log 2+S_p(\infty)>0.
 \]
\par
 Thus, if
\[
\begin{aligned}
t>C_7\varepsilon^{-2p(p-1)/\gamma(p, n)}
\end{aligned}
\]
with
\[
C_7:=\Big(\frac{n(n+1)e^{C_6+1}}{C_4}\Big)^{2(p-1)/\gamma(p, n)},
\]
we then get $J(t)>1$,
and this in turn gives that $F_0(t)\rightarrow \infty$ by letting $j\rightarrow \infty$ in \eqref{27}. Therefore we get the desired upper bound,
\[
\begin{aligned}
T\leq C_7\varepsilon^{-2p(p-1)/\gamma(p, n)},
\end{aligned}
\]
and hence we finish the proof of Theorem \ref{thm1}.


\section{Proof for Theorem \ref{thm2} and Theorem \ref{thm3}}

Due to \eqref{7a}, (\ref{9}) and (\ref{12}), integrating \eqref{11} over $[0, t]$ yields
\[
\begin{aligned}
F_0'(t)\geq m(t)F_0'(t)\geq C_8\e,
\end{aligned}
\]
where
\[
C_8:=m(0)\int_{\R^n}g(x)dx.
\]
The above inequality implies that
\begin{equation}
\label{81}
F_0(t)\geq C_9\e(1+t)\quad\mbox{for}\quad t\geq 0,
\end{equation}
where
\[
C_9:=\min\left\{C_8,\int_{\R^n}f(x)dx\right\}.
\]
\par
First we prove Theorem \ref{thm2} for $n=2$.
Due to the assumption on $g(x)$, we note that
\[
\int_{\R^2}g(x)dx>0
\]
which yields $C_8,C_9>0$.
By (\ref{19a}) and \eqref{81}, we have
\begin{equation}
\label{81a}
\int_{\R^2}|u(x, t)|^pdx\geq C_{10}\e^p(1+t)^{2-p},
\end{equation}
with $C_{10}:=C_2C_9^p>0$.
Plugging \eqref{81a} into \eqref{12} and integrating it over $[0, t]$ we come to
\begin{equation}
\label{83}
\begin{aligned}
F_0(t)&\geq C_{11}\e^p\int_0^tdr\int_0^r(1+s)^{2-p}ds\\
&\geq C_{11}\e^p(1+t)^{1-p}\int_0^tdr\int_0^rsds\\
&=\frac{C_{11}}{6}\e^p(1+t)^{1-p}t^3
\quad\mbox{for}\quad t\ge0
\end{aligned}
\end{equation}
with $C_{11}:=C_{10}m(0)>0$. Noting that the above inequality improves the lower bound of \eqref{20} for $n=2$ and $1<p<2$,
and this is the key point to prove Theorem \ref{thm2}.

As in section 4, we define our iteration sequence,
$\{\widetilde{D}_j\}, \{\widetilde{a}_j\}, \{\widetilde{b}_j\}$, as
\begin{equation}\label{85}
\begin{aligned}
F_0(t)\geq \widetilde{D}_j(1+t)^{-\widetilde{a}_j}t^{\widetilde{b}_j}
\quad\mbox{for}\quad t\geq 0, \ j=1, 2, 3\cdots
\end{aligned}
\end{equation}
with positive constants, $\widetilde{D}_j, \widetilde{a}_j, \widetilde{b}_j$, and
\begin{equation}\nonumber
\begin{aligned}
\widetilde{D}_1=\frac{C_{11}}{6}\e^p,~~\widetilde{a}_1=p-1,~~\widetilde{b}_1=3.
\end{aligned}
\end{equation}
Combining (\ref{19b}) and \eqref{85}, we have
\[
\begin{aligned}
F_0(t)&\geq C_3\widetilde{D}_j^p\int_0^tdr\int_0^r(1+s)^{2-2p-p\widetilde{a}_j}s^{p\widetilde{b}_j}ds\\
&\geq \frac{C_3\widetilde{D}_j^p}{(p\widetilde{b}_j+2)^2}(1+t)^{2-2p-p\widetilde{a}_j}t^{p\widetilde{b}_j+2}
\quad\mbox{for}\quad t\geq 0.\\
\end{aligned}
\]
So the sequences satisfy
\begin{equation}\nonumber
\begin{aligned}
\widetilde{a}_{j+1}&=-p\widetilde{a}_j-2(p-1),\\
\widetilde{b}_{j+1}&=p\widetilde{b}_j+2,\\
\widetilde{D}_{j+1}&\geq \frac{C_3\widetilde{D}_j^p}{(p\widetilde{b}_j+2)^2}\geq \frac{C_{12}\widetilde{D}_j^p}{p^{2j}},
\end{aligned}
\end{equation}
where $C_{12}:=C_{11}/\{3+2/(p-1)\}^2>0$, which means that
\begin{equation}\nonumber
\begin{aligned}
\widetilde{a}_{j}&=p^{j-1}(p+1)-2,\\
\widetilde{b}_{j}&=\big(3+\frac{2}{p-1}\big)p^{j-1}-\frac{2}{p-1},\\
\log\widetilde{D}_{j}&\geq p^{j-1}\log\widetilde{D}_1-\sum_{k=1}^{j-1}\frac{2k\log p-\log C_{12}}{p^k}.
\end{aligned}
\end{equation}
Then, as in Section 4, we have
\[
\begin{aligned}
F_0(t)&\geq\widetilde{D}_j(1+t)^{-p^{j-1}(p+1)+2}
t^{p^{j-1}\{3+2/(p-1)\}-2/(p-1)}\\
&\geq (1+t)^2t^{-2/(p-1)}\exp\big(p^{j-1}\widetilde{J}(t)\big),
\end{aligned}
\]
where
\[
\widetilde{J}(t):=-(p+1)\log(1+t)+\left(3+\frac{2}{p-1}\right)\log t+\log \widetilde{D}_1-\widetilde{S}_p(\infty)
\]
and
\[
\widetilde{S}_p(\infty):=\sum_{k=1}^{\infty}\frac{2k\log p-\log C_{12}}{p^k}.
\]
Estimating $\wt{J}(t)$ as
\[
\begin{aligned}
\wt{J}(t)
&\geq -(p+1)\log(2t)+\left(3+\frac{2}{p-1}\right)\log t+\log \widetilde{D}_1-\widetilde{S}_p(\infty)\\
&=\left(-\frac{p(p-3)}{p-1}\right)\log t+\log \widetilde{D}_1-\widetilde{S}_p(\infty)-(p+1)\log 2,
\end{aligned}
\]
we obtain that 
\[
\wt{J}(t)\ge\log\Big(t^{-p(p-3)/(p-1)}\widetilde{D}_1\Big)-C_{13}
\quad\mbox{for}\quad t\geq 1,
\]
where $C_{13}:=\widetilde{S}_p(\infty)+(p+1)\log 2>0$.
By the definition of $\widetilde{D}_1$, we then get the lifespan estimate in Theorem \ref{thm2}
by the same way as that in section 4.

\par
Next we prove Thorem \ref{thm3} for $n=1$.
The proof can be shown along the same way as that of Theorem \ref{thm2} for $n=2$.
First we note that (\ref{81}) is also available in this case.
Then the first iteration in \eqref{83} for $n=2$ becomes
\[
F_0(t)\geq C_{11}\e^p\int_0^tdr\int_0^r(1+s)ds
\geq \frac{C_{11}}{6}\e^p t^3
\quad\mbox{for}\quad t\geq 0.
\]
We then may assume the iteration sequences, $\{\o{D}_j\},\{\o{a}_j\},\{\o{b}_j\}$, as
\[
F_0(t)>\o{D}_j(1+t)^{-\o{a}_j}t^{\o{b}_j}\quad\mbox{for}\ t\ge0
\quad(j=1,2,3\cdots)
\]
with non-negative constants, $\o{D}_j,\o{a}_j,\o{b}_j$, and
\begin{equation}\nonumber
\o{D}_1=\frac{C_{11}}{6}\e^p,
\quad \o{a}_1=0,
\quad \o{b}_1=3.
\end{equation}
The left steps to get the lifespan in Theorem \ref{thm3} are exactly the same as that of Theorem \ref{thm2}.

\section*{Acknowledgment}
\par\quad
The first author is partially supported by NSFC(11501273, 11771359), high level talent project
of Lishui City(2016RC25), the Scientific Research Foundation of the First-Class Discipline of Zhejiang Province(B)(201601), the key laboratory of Zhejiang Province(2016E10007).
The second author is partially supported by the Grant-in-Aid for Scientific Research(C)
(No.15K04964),
Japan Society for the Promotion of Science,
and Special Research Expenses in FY2017, General Topics(No.B21),
Future University Hakodate.
\par
Finally, the authors are grateful to Prof.M.Ikeda
for his pointing out the mistake in the proof for the critical case
in  our first draft at arXiv:1707.09583.


\bibliographystyle{plain}

\end{document}